\documentclass[12 pt]{amsart}
\usepackage{amscd,amssymb,amsmath,amsthm}
\newtheorem{Proposition}{Proposition}

\newtheorem{Lemma}{Lemma}
\newtheorem{Theorem}{Theorem}
\newtheorem{Corollary}{Corollary}

\newcommand{\proj}{\mathbb{P}}

\newcommand{\rarr}{\rightarrow}
\newcommand{\oh}{{\mathcal{O}}}
\newcommand{\com}{\mathbb{C}}
\newcommand{\Q}{\mathbb{Q}}

\newcommand{\Lan}{\Big \langle}
\newcommand{\Ran}{\Big \rangle}

\begin{document}

\title{New calculations in Gromov-Witten theory}
\author{D.~Maulik and R.~Pandharipande}
\dedicatory{Dedicated to Fedor Bogomolov on the occasion of his 60th birthday}
\date{January 2006}
\maketitle

\setcounter{section}{-1}
\section{Introduction}
\subsection{Overview}
Let $X$ be a nonsingular projective variety over $\com$.
The stack $\overline{M}_{g,n}(X,\beta)$ parameterizes
stable maps 
$$f: C \rightarrow X$$
from genus $g$, $n$-pointed curves to $X$ representing the class
$\beta\in H_2(X,{\mathbb{Z}})$.
The Gromov-Witten invariants of $X$ are defined by 
integration
against
the virtual class 
of the moduli
space,
\begin{equation}
\label{xty}
\Big\langle \tau_{a_1}(\gamma_1)\cdots \tau_{a_n}(\gamma_n) 
\Big\rangle^X_{g,\beta} = \int_{[\overline{M}_{g,n}(X,\beta)]^{vir}}
\prod_{i=1}^n \psi_i^{a_i} \cup \text{ev}_i^*(\gamma_i),
\end{equation}
where 
$\psi_i$ is the Chern class of the $i^{th}$ cotangent line,
$$\text{ev}_i: \overline{M}_{g,n}(X,\beta) \rightarrow X$$
is the $i^{th}$ evaluation map, and
$\gamma_i\in H^*(X,\mathbb{Q})$.

Over the past 10 years, there has been great deal of success in calculating
Gromov-Witten invariants if {either} the domain genus is 0 {or}
the target $X$ carries a strong torus action --- a torus action
with finitely many  0 and 1 dimensional orbits.

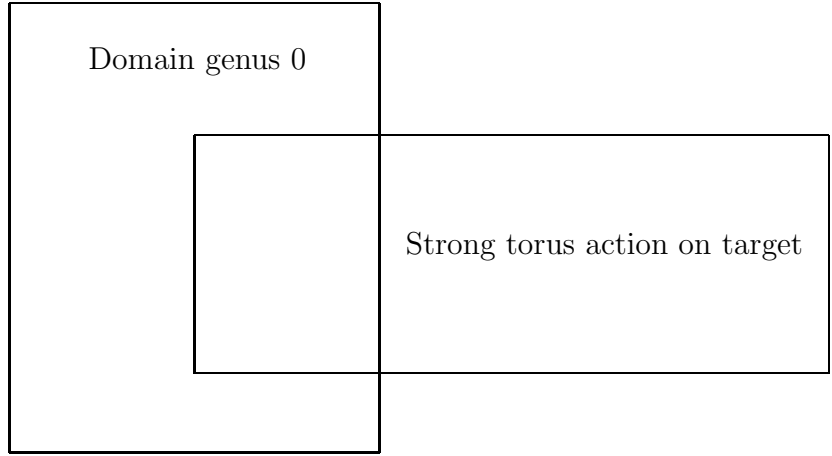
\begin{picture}(200,200)(-20,-5)
\put(10,10){\line(1,0){140}}
\put(80,40){\line(1,0){240}}
\put(80,130){\line(1,0){240}}
\put(10,180){\line(1,0){140}}
\put(10,10){\line(0,1){170}}
\put(80,40){\line(0,1){90}}
\put(320,40){\line(0,1){90}}
\put(150,10){\line(0,1){170}}
\put(40,155){Domain genus 0}
\put(160,85){Strong torus action on target}
\thicklines
\end{picture}

\noindent The goal of our paper is to use the methods of \cite{mp} 
for Gromov-Witten calculations outside of the above regions.

We study the higher genus Gromov-Witten invariants
of two target geometries: surfaces of general 
type and compact Calabi-Yau 3-folds. In the surface case,
the calculations suggest
exact solutions. For Calabi-Yau 
3-folds, our results provide
the first mathematical encounter with the holomorphic
anomaly equation for topological strings.

\subsection{Genus 0}
There are several mathematical
approaches
to Gromov-Witten theory in genus 0.
A wide class
of target varieties has been successfully studied
via WDVV-equations, Frobenius structures, and quantum Lefschetz
formulas. Quantum Schubert calculus is well developed for
classical homogeneous spaces. Mirror symmetry relations between
genus 0 invariants of Calabi-Yau hypersurfaces and Picard-Fuchs
systems have been proven in many cases. 
The
bibliography of \cite{ck} is a good source for
genus 0 references before 2000. 

An example of recent progress in genus 0 is the calculation
of the
quantum cohomology of the Hilbert scheme of points
of $\com^2$ \cite{qchs}.

Genus 0 computation for {\em orbifold} targets, however, is a largely open
subject \cite{bgr,bgrp,cr}.

\subsection{Strong torus actions}
The situation in higher genus is substantially more difficult. If $X$ carries a
strong torus action, then the Gromov-Witten theory of $X$
is reduced to Hodge integrals on $\overline{M}_{g,n}$
via localization of the virtual class \cite{GP}.
Hodge integrals can be reduced to descendent integrals by Gromov-Witten 
operators based on a Grothendieck-Riemann-Roch calculation
of Mumford \cite{FP}. 
A formalism for expressing the higher genus invariants of $X$ in terms
of genus 0 Frobenius structures via \cite{FP,GP} has
been developed by Givental \cite{giv,yp}. An outcome, for example, 
is Givental's proof of the Virasoro constraints
for the Gromov-Witten theory of $\mathbb{P}^n$.

If the dimension of $X$ is 3, calculations in higher genus 
can be made much more effective. The required
Hodge integrals
have been treated with 
increasing sophistication. 
The culmination has been the topological vertex \cite{topv} 
in the local Calabi-Yau toric case and the
equivariant vertex \cite{mnop1,mnop2} for arbitrary 3-folds with strong torus
actions.
Neither vertex evaluation of 3-fold Gromov-Witten theory 
has yet a complete
mathematical proof, see \cite{topv,llz,mnop1,mnop2,cry}.

\subsection{New directions}
If we leave the genus 0 or strong torus action realm, very
few calculations have been completed. 

The
Gromov-Witten invariants of all 1-dimensional targets 
have been fully determined in \cite{OP1,OP2,OP3}.
The method uses a mix of localization, degeneration,
and exact evaluations. Similarly, the local theory of curves
in 3-folds has been solved in \cite{bp}.
The results of \cite{mp} present
a topological view of Gromov-Witten calculations in all
dimensions: the familiar cutting and pasting strategies
in the topological category
are shown to yield effective Gromov-Witten
techniques.

After a discussion of our topological view in Section \ref{one},
applications to branched double covers of $\mathbb{P}^2$  and 
the Enriques Calabi-Yau 3-fold are
presented in Sections \ref{two} and \ref{three}.

\subsection{Acknowledgments}
Conversations with J. Bryan, T. Parker, and M. Usher played an
important role in our surface calculations. Many of these
occurred during a workshop on holomorphic curves at the
Institute for Advanced Study in Princeton in June of 2005.

Our work on the Enriques Calabi-Yau was largely motivated by a
series of lectures by A. Klemm and M. Mari\~no  on heterotic
duality and the holomorphic anomaly equation for the
Enriques geometry. The lectures (and most of our calculations)
took place at a workshop on algebraic geometry and topological
strings at the Instituto Superior T\'ecnico in Lisbon in the fall
of 2005. 

We thank I. Dolgachev for his immediate answers to all of our
questions about the classical geometry of the Enriques surface.

D.~M. was partially supported by an NSF graduate fellowship.
R.~P. was partially supported by the NSF and the Packard
foundation.
\pagebreak
\section{Topological view}
\label{one}
\subsection{Mayer-Vietoris}
\label{oneone}
Let $X$ be a nonsingular projective variety. We would like to
study Gromov-Witten theory by decomposing $X$ into simpler
pieces.

Degeneration may be viewed as an algebraic version of cutting
and pasting. 
Let
$$\epsilon:{\mathcal X} \rarr \Delta$$
be a flat family over a disk $\Delta \subset \com$ at the origin
satisfying:
\begin{enumerate}
\item[(i)] $\mathcal {X}$ is nonsingular,
\item[(ii)] $\epsilon$ is smooth over the punctured disk $\Delta^*=\Delta\setminus \{0\}$,
\item[(iii)] $\epsilon^{-1}(1)\stackrel{\sim}{=} X$,
\item[(iv)] $\epsilon^{-1}(0)=X_1\cup_Y X_2$ is a normal crossings divisor in ${\mathcal X}$.
\end{enumerate}
The family $\epsilon$ defines a canonical map
$$H^*(X_1\cup_Y X_2,\Q) \rarr H^*(X,\Q)$$
with image defined to be the {\em nonvanishing} cohomology of $X$.

Let $\mathsf{GW}(X)$ denote the 
descendent Gromov-Witten theory of the target $X$ --- the complete
set of integrals \eqref{xty} --- and let $\mathsf{GW}(X)_\epsilon$
denote the Gromov-Witten theory of descendents of the
nonvanishing cohomology.
A Mayer-Vietoris result in Gromov-Witten theory is proven for the
family ${\mathcal{X}}$
in \cite{mp}. 

\begin{Theorem}\label{mv}
$\mathsf{GW}(X)_\epsilon$
can be uniquely and effectively
reconstructed from 
$$\mathsf{GW}(X_1), \ \mathsf{GW}(X_2), \ \mathsf{GW}(Y),$$
 and the
restriction maps
$$H^*(X_1,\Q) \rarr H^*(Y,\Q)\leftarrow  H^*(X_2,\Q).$$ 
\end{Theorem}

A degeneration $\epsilon:{\mathcal X} \rightarrow \bigtriangleup$
satisfying conditions (i)-(iv) above will be called {\em good}.
Theorem \ref{mv} requires a good degeneration.

\subsection{Surfaces branched over $\mathbb{P}^2$} \label{on2}
Let $S$ be a nonsingular projective surface which admits
a good degeneration to $S_1 \cup_C S_2$.
Since $\mathsf{GW}(C)$  has
been completely determined in \cite{OP3}, Theorem \ref{mv}
is particularly applicable.

Surfaces constructed as branched double covers of $\mathbb{P}^2$
provide a basic class of  
examples. 
Let
$$S_{2n} \rightarrow \mathbb{P}^2$$
be a double cover branched along a nonsingular curve $B$ of
degree $2n$.
If $n\geq 4$, then $S_{2n}$ is of general type.

Let $C$ be a nonsingular plane curve of degree $n$ generic
with respect to $B$.
By
degenerating the branch curve $B$ to the square of $C$, 
we can construct
a 1-parameter family ${\mathcal{F}}$ of surfaces:
$${\mathcal{F}} \rightarrow {\mathbb{P}^2} \times \com$$
is the double cover along
$$(tB-C^2)\subset {\mathbb{P}^2} \times \com$$
where $t$ is the parameter on $\com$.
The total space ${\mathcal{F}}$ has double point singularities
above the $2n^2$ points of
$B \cap C$. Taking the small resolution, we obtain a
nonsingular space
$$\epsilon:{\mathcal S} \rightarrow \com$$
which provides a degeneration of $S_{2n}$ to
$$\tilde{\mathbb{P}}^2 \cup_C \mathbb{P}^2$$
where $\tilde{\mathbb{P}}^2$ is the blow-up of
$\mathbb{P}^2$ along $B\cap C$.

The cohomology of $S_{2n}$ pulled-back from $\mathbb{P}^2$
is certainly nonvanishing for $\epsilon$.
Since the descendent theories $\mathsf{GW}(\tilde{\mathbb{P}}^2)$
and $\mathsf{GW}({\mathbb{P}}^2)$ are accessible via
various methods (localization, Virasoro, Frobenius structures),
Theorem \ref{mv} provides an effective approach to 
$\mathsf{GW}(S_{2n})_\epsilon$.

\subsection{The Enriques surface} \label{13s}
Let $\sigma_{K3}$ be a fixed point free involution of a K3 surface.
By definition, the quotient
$$X = K3/\langle\sigma_{K3}\rangle$$ 
is an {\em Enriques surface}.
Alternatively, Enriques surfaces arise as elliptic fibrations
\begin{equation}\label{dfh1}
X \rightarrow \mathbb{P}^1
\end{equation}
with 12 singular fibers and 2 double fibers.

Let $E\times \mathbb{P}^1$ be the product of an elliptic curve 
and a projective line. Let $\text{t}_E$ denote translation on $E$
by a 2-torsion point, and let $\text{inv}_{\mathbb{P}^1}$ denote
an involution on the projective line.
Then,
$$\tau: E \times \mathbb{P}^1 \rightarrow E \times \mathbb{P}^1$$
acting by $(\text{t}_E,\text{inv}_{\mathbb{P}^1})$
is a fixed point free involution. Let
$$X_1 = \left( E \times \mathbb{P}^1 \right)/\langle \tau \rangle$$
denote the quotient.
By projecting left,  
$$X_1 \rightarrow E/\langle \text{t}_E \rangle$$
is a projective bundle over the elliptic curve $E/\langle \text{t}_E \rangle$.
Hence $\mathsf{GW}(X_1)$ is determined by localization and
\cite{OP3}. By projecting right,
$$X_1 \rightarrow \mathbb{P}^1/\langle \text{inv}_{\mathbb{P}^1}\rangle$$
is an elliptic fibration with no singular fibers and 2 double
fibers.

Let $E(1)$ be the {\em rational elliptic surface} isomorphic
to the blow-up of $\mathbb{P}^2$ in 9 points. The surface $E(1)$
admits an elliptic fibration
$$E(1) \rightarrow \mathbb{P}^1$$
with 12 singular fibers and no double fibers.

By degenerating the elliptic fibration \eqref{dfh1}, we find a good
degeneration of the Enriques surface $X$ to
$$X_1 \cup _E E(1)$$
where the intersection $E$ is a common elliptic fiber.
Since all the cohomology of $X$ is nonvanishing for the
degeneration, Theorem
\ref{mv} provides an effective approach to the Gromov-Witten
theory of the Enriques surface.

Since the $K3$ surface is holomorphic symplectic,
$\mathsf{GW}(K3)$ is essentially 
trivial\footnote{The $K3$ surface has a rich {\em modified}
Gromov-Witten theory. The investigation of the modified Gromov-Witten
theory of the $K3$ surface
 should also be considered to be outside of the genus 0 and
toric realm. Calculations in primitive and twice primitive
classes can be found in \cite{brl,ll}. At present, the
modified virtual class is not amenable to the techniques
discussed here.}
with nonvanishing invariants only for constant maps
in genus 0 and 1.
The Enriques surface, however, will be seen to have a very
rich Gromov-Witten theory.

\subsection{The Enriques Calabi-Yau}\label{et45}
Let $\sigma$ act freely on the product $K3 \times E$ by
an Enriques involution $\sigma_{K3}$ on the $K3$ and by -1
on the elliptic curve.
By definition, the quotient
\begin{equation*}
Q = \left( 
K3 \times E\right) / \langle \sigma 
\rangle
\end{equation*}
is an {\em Enriques Calabi-Yau 3-fold}.
Since $K3\times E$ carries a holomorphic 3-form invariant
under $\sigma$,
$$K_Q =0.$$
By projection on the right,
\begin{equation}\label{ft4}
Q \rightarrow E/\langle -1 \rangle = \mathbb{P}^1
\end{equation}
is a $K3$ fibration with 4 double Enriques fibers.

Let $\tau$ act freely on the product $K3 \times {\mathbb{P}}^1$ by
$(\sigma_{K3},\text{inv}_{\mathbb{P}^1})$.
Let
$$R = \left( 
K3 \times \mathbb{P}^1 \right) / \langle \tau 
\rangle$$
denote the quotient.
By projecting left,  
$$R \rightarrow K3/\langle \sigma_{K3} \rangle=X$$
is a projective bundle over the Enriques surface $X$.
Hence $\mathsf{GW}(R)$ is determined by localization and
$\mathsf{GW}(X)$. By projecting right,
$$R \rightarrow \mathbb{P}^1/\langle \text{inv}_{\mathbb{P}^1}\rangle$$
is a K3 fibration with  2 double Enriques
fibers.

By degenerating the $K3$ fibration \eqref{ft4}, we find a
good degeneration of the Enriques Calabi-Yau $Q$ to
$$R \cup _{K3} R$$
where the intersection $K3$ is a common fiber.
Since all the cohomology of $Q$ is nonvanishing, Theorem
\ref{mv} provides an effective approach to the Gromov-Witten
theory of the Enriques Calabi-Yau.

\subsection{Absolute/Relative}
Relative Gromov-Witten theory plays an essential role in
Theorem \ref{mv}. We review standard notation for relative
invariants.

Let $(V,W)$ be a nonsingular projective variety $V$ containing
a nonsingular divisor $W$.
Let $\beta\in H_2(V,{\mathbb Z})$ be a curve class satisfying
$$\int_\beta[W]\geq 0.$$
Let $\stackrel{\rightarrow}{\mu}$ be an ordered partition,
$$\sum_j \mu_j = \int_\beta[W],$$
with positive parts.
The moduli space 
$\overline{M}_{g,n}(V/W,
\beta, \stackrel{\rightarrow}{\mu})$ parameterizes
stable relative 
maps from genus $g$, $n$-pointed curves to $V$ of class $\beta$
with 
multiplicities along $W$ determined by $\stackrel{\rightarrow}{\mu}$.

The relative 
conditions in the theory correspond to partitions {\em weighted} by
the cohomology of $W$.
Let  $\delta_1, \ldots, \delta_{m_W}$ be a basis of $H^*(W,{\mathbb Q})$.
A cohomology weighted partition ${\nu}$
consists of an {\em unordered} set of pairs,
$$\left\{ (\nu_1, \delta_{s_1}), \ldots, 
(\nu_{\ell(\nu)}, \delta_{s_{\ell(\nu)}}) 
\right\}, $$
where $\sum_j \nu_j$ is an {\em unordered} partition of $\int_\beta [W]$.
The automorphism group, $\text{Aut}(\mathbf{\nu})$, consists of 
permutation symmetries of ${\mathbf \nu}$.

The {\em standard} order on the parts of $\nu$ is
$$(\nu_i,\delta_{s_i})> (\nu_{i'},\delta_{s_{i'}})$$
if $\nu_i>\nu_{i'}$ or if $\nu_i=\nu_{i'}$ and $s_i>s_{i'}$.
Let 
$\stackrel{\rarr}{\nu}$ denote the partition $(\nu_1,\ldots, \nu_{\ell(\nu)})$
obtained from the standard order.

The descendent Gromov-Witten invariants of 
the pair $(V,W)$ are defined by integration against the
virtual class of the moduli of maps. 
Let $\gamma_1,\ldots, \gamma_{m_V}$ be a basis of $H^*(V,\Q)$, and
let
\begin{multline}\label{f32}
\left. \Lan  \tau_{k_1}(\gamma_{l_1}) \cdots
\tau_{k_n}(\gamma_{l_n})\ \right| {\mathbf \nu} 
\Ran^{V/W}_{g,\beta} = \\
 \frac{1}{|\text{Aut}(\nu)|}
\int_{[\overline{M}_{g,n}(V/W,\beta,\stackrel{\rarr}{\nu})]^{vir}} 
\prod_{i=1}^n \psi_i^{k_i}\text{ev}_i^*(\gamma_{l_i}) \cup 
\prod_{j=1}^{\ell(\nu)} \text{ev}^*_j(\delta_{s_j}).
\end{multline}
Here, the second evaluations, 
$$\text{ev}_j: \overline{M}_{g,n}(V/W,\beta,\stackrel{\rarr}\nu) \rarr W.$$
are determined by the relative points.
The brackets \eqref{f32} denote integration over the
moduli of maps with {\em connected} domains.

Gromov-Witten invariants are defined (up to sign) for {\em unordered}
weighted partitions ${\nu}$. To fix the sign, the integrand
on the right side requires an ordering. The ordering is corrected
by the automorphism prefactor.

Given a good degeneration of $X$ to $X_1\cup_Y X_2$,
the degeneration formula \cite{EGH,IP,LR,L} expresses
$\mathsf{GW}(X)_\epsilon$ in terms of the relative
Gromov-Witten theories of the pairs
$(X_1,Y)$ and $(X_2,Y)$. Hence, Theorem \ref{mv} is 
a consequence of the following result of \cite{mp}.

\begin{Theorem}
\label{wd}
The relative Gromov-Witten 
theory of the pair $(V,W)$ can be uniquely and effectively
reconstructed from $\mathsf{GW}(V)$, $\mathsf{GW}(W)$, 
and the restriction map
$H^*(V,\Q) \rarr H^*(W,\Q)$.
\end{Theorem}

\section{Surfaces of general type}\label{two} 
\subsection{Seiberg-Witten theory}
Let $S$ be a nonsingular, projective, minimal surface of
general type with $p_g(S)>0$. Let $K_S\in H_2(S,{\mathbb{Z}})$ be the canonical
class of $S$, and let
$$g_K = K_S^2+1$$
be the adjunction genus in the canonical class.
The moduli space $\overline{M}_{g_K}(S,K_S)$ has virtual
dimension 0. 
Taubes has obtained the
evaluation
\begin{equation}
\label{ct5}
\langle 1 \rangle_{g_K,K_S}^S = (-1)^{\chi(\mathcal{O}_S)}.
\end{equation}
by a connection to Seiberg-Witten theory \cite{iprt,morgan,t1,t2,t3,t4}.
Here, $\chi(\mathcal{O}_S)$ denotes the holomorphic Euler
characteristic.

We will call a nonsingular, irreducible, canonical divisor $C\subset S$ 
a {\em Taubes curve}.
If a Taubes curve exists,
formula \eqref{ct5} has a simple geometric
interpretation. By adjunction, the normal bundle
of $C$ in $S$ is a square root of the canonical
bundle of $C$,
\begin{eqnarray}
\oh_C(2C) & = & K_S(C)|_C \label{sd4}\\
 & = & K_C.\nonumber 
\end{eqnarray}
The normal bundle is thus a {\em theta characteristic}
of $C$ and has a deformation invariant parity equal to
$h^0(C,\oh_C(C)) \mod 2$.
The {\em sign} of the Taubes curve is defined by the
parity of the normal bundle,
$$\sigma(C)= (-1)^{h^0(C,\oh_C(C))}.$$

\begin{Lemma}\label{d45} If $C\subset S$ is a Taubes curve,
then $\langle 1 \rangle_{g_K,K_S}^S$ equals the sign of $C$.
\end{Lemma}

\begin{proof} The cohomology sequence associated to the
short exact sequence
$$0 \rarr \oh_S \rarr \oh_S(C) \rarr \oh_C(C) \rarr 0$$
starts as
\begin{eqnarray*}
0 &\rarr &H^0(S,\oh_S) \rarr H^0(S,\oh_S(C)) \rarr H^0(C,\oh_C(C))
\rarr \\
& & H^1(S,\oh_S) \stackrel{\phi}{\rarr} H^1(S,\oh_S(C)) \rarr \ldots \, .
\end{eqnarray*}
Since $C$ is a Taubes curve,
$$h^0(S,\oh_S(C))=h^0(S,K_S)=h^2(S,\oh_S).$$
Hence, the parity of the normal bundle equals the parity of
$$\chi(\oh_S)+ \text{dim}(\text{Im}(\phi)).$$
The image $\text{Im}(\phi)$ has dimension equal to the rank of
the skew-symmetric form on $H^1(S,\oh_S)$ defined
by
$$H^1(S,\oh_S) \times H^1(S,\oh_S) \stackrel{(\phi,id)}{\rarr}
H^1(S,K_S) \times H^1(S,\oh_S) \stackrel{\cup}{\rarr}
H^2(S,K_S) \stackrel{\sim}{=}\com.$$ 
The rank of a skew-symmetric form is even.
\end{proof}

Lemma \ref{d45} is well-known in various forms. Our purpose
is to emphasize the connection to the theta
characteristic.

\subsection{Vanishing and universality}
The Gromov-Witten
invariants of $S$ of {\em Severi} type  are
$$\Big \langle \tau^n_0([p_S]) 
\Big \rangle_{g,\beta\neq 0}^S,$$
where $[p_S]\in H^4(S,\mathbb{Z})$ is the point class.
By results of Taubes,
Gromov-Witten invariants of
Severi type vanish in the {\em adjunction genus}
$$2g_\beta-2 = (K_S+\beta)\cdot \beta$$ 
 if 
$\beta \neq K_S$.
The universal evaluation
$$\langle 1 \rangle_{g_K,K_S}^S = (-1)^{\chi(\mathcal{O}_S)}$$
is an example of nonvanishing.

Let $S$ be a minimal surface of general type.
If a
Taubes curve $C\subset S$ exists,
a stronger vanishing result is proven by
J. Lee and T. Parker 
in \cite{park}. The methods of \cite{park} do not
use Seiberg-Witten theory.

\vspace{10pt}
\noindent{\bf I. Vanishing:} {\em 
The Gromov-Witten invariants
$$\Big \langle \tau_{\alpha_1}(\gamma_1) \ldots \tau_{\alpha_n}(\gamma_n) 
\Big \rangle_{g,\beta\neq 0}^S$$
vanish if either $\beta \notin \mathbb{Z}K_S$ or there exists
an insertion satisfying
$$\gamma_i \in H^{\geq 3}(S,\mathbb{Q}).$$}

Constant maps are avoided in the vanishing statement. A
complete discussion of $\beta=0$ invariants of a surface can
be found in \cite{getp}.

\vspace{9pt}
\noindent{\bf II. Universality:} {\em 
$\mathsf{GW}(S)$ is uniquely determined
by the sign of the Taubes curve and the restriction map
$H^*(S,\mathbb{Q}) \rightarrow H^*(C,\mathbb{Q})$.}
\vspace{9pt}

The universality II is slightly stronger than the result
of \cite{park} where the precise spin structure \eqref{sd4}
is required. Since the sign is the only deformation
invariant of a spin structure, we conjecture
the stronger universality II.

\subsection{Local theory of surfaces}
We will assume the existence of a
Taubes curve $C\subset S$. A sharper
universality statement can be made for 
 descendents of the even cohomology 
$$H^0(S,\mathbb{Q}) \oplus H^2(S,\mathbb{Q})\subset H^*(S,\mathbb{Q}).$$
of $S$. By vanishing, descendents of $H^4(S,\mathbb{Q})$
need not be considered.

\vspace{9pt}
\noindent{\bf II$'$. Universality:} 
For $d>0$,
\begin{multline*}
\Big \langle \prod_{i=1}^n
\tau_{\alpha_i}(D_i) \prod_{j=1}^m 
\tau_{\tilde{\alpha}_j}(1) 
\Big \rangle_{g,dK_S}^S= \\
d^n\prod_{i=1}^n ({K_S}\cdot D_i) \ \cdot \ 
\Big \langle \prod_{i=1}^n
\tau_{\alpha_i}([p_C]) \prod_{j=1}^m 
\tau_{\tilde{\alpha}_j}(1) 
\Big \rangle_{g,d}^{C,\sigma(C)},
\end{multline*}
where $D_i \in H^2(S,\mathbb{Q})$ and
$[p_C]\in H^2(C,\mathbb{Q})$ is the point class.
\vspace{9pt}

The brackets $\langle,\rangle_{g,d}^{C,\pm}$
refer to a {\em local} Gromov-Witten theory of curves
in surfaces for which no explicit algebraic construction yet exists.
The universality equation II$'$ may be taken to
conjecturally
{\em define} the local theory on the right. 

\subsection{Conjectures}
We conjecture evaluations of the local theory of curves in surfaces
\begin{equation} \label{fgq1}
\Big \langle \prod_{i=1}^n
\tau_{\alpha_i}([p_C]) \prod_{j=1}^m 
\tau_{\tilde{\alpha}_j}(1) 
\Big \rangle_{g,d}^{C,\pm}
\end{equation}
 for degrees $d=1,2$. 
The dimension constraint for the local theory \eqref{fgq1}
is
$$g-1-d(g_C-1)+m= \sum_{i=1}^n \alpha_i + \sum_{j=1}^m \tilde{\alpha}_j.$$
By the Virasoro conjecture for the Gromov-Witten theory
of surfaces,\footnote{
The Virasoro conjecture is open, but
partial results for the local theory of surfaces
have been obtained by A. Gholampour.} the insertions $\tau_{\tilde{\alpha}}(1)$
can be removed by universal relations. 

We will 
restrict our attention to the insertions $\tau_{\alpha}([p_C])$.
In degree 1, we conjecture
\begin{equation}\label{fg3}
\Big \langle \prod_{i=1}^n
\tau_{\alpha_i}([p_C]) 
\Big \rangle_{g,1}^{\bullet\ C,\pm}=
\pm \prod_{i=1}^n \frac{\alpha_i !}{(2\alpha_i+1)!} (-2)^{-\alpha_i}.
\end{equation}
In degree 2, we conjecture
\begin{equation}\label{fg4}
\Big \langle \prod_{i=1}^n
\tau_{\alpha_i}([p_C])  
\Big \rangle_{g,2}^{\bullet\ C,\pm}=
\pm \ 2^{g_C+n-1} \prod_{i=1}^n \frac{\alpha_i !}{(2\alpha_i+1)!} 
(-2)^{\alpha_i}.
\end{equation}
The superscripted bullet here denotes
Gromov-Witten invariants with possibly disconnected domains
{\em with no degree 0 components}.
The heuristic origins of formulas \eqref{fg3} and \eqref{fg4} 
are based on relationships
to exact evaluations in the Gromov-Witten theory of $\mathbb{P}^1$.
A discussion will be presented in \cite{dm}.

Our main evidence for the degree 1 and 2 formulas is obtained
from calculations in the Gromov-Witten theory of $S_{2n}$
based on the strategy of Section \ref{on2}. The computations are
rather labor intensive --- very similar to the quintic 
surface computations of \cite{mp} --- and will be omitted
here. 

For readers seeking worked examples of the method of \cite{mp}, the 3-fold
calculations  of Section
\ref{three} will be explained in detail.

\subsection{Speculations}
The full evaluation of the local theory 
 $$\Big\langle\ ,\ \Big\rangle_{g,d}^{C,\pm}$$
of curves in
surfaces 
is not yet known.
The Gromov-Witten theory
of curves is related to gauge theory for the symmetric group
via the Gromov-Witten/Hurwitz correspondence of \cite{OP1}.
The local Gromov-Witten
theory of curves in 3-folds is related to $U(1)$-gauge
theory via the Gromov-Witten/Donaldson-Thomas correspondence
\cite{bp,opnew}.
The local Gromov-Witten theory of curves in surfaces sits
in between. 
 We expect the local theory
of curves in surfaces to be related to a descendent Donaldson 
theory of sheaves on surfaces.

\section{The Enriques Calabi-Yau 3-fold}
\label{three}

\subsection{Fiber classes}
\subsubsection{2-torsion}
Let $X$ be the Enriques surface.
The middle homology of $X$ has a 2-torsion factor. 
The torsion-free quotient,
$$H_2(X,\mathbb{Z})' \stackrel{\sim}{=}
H_2(X,\mathbb{Z})/\mathbb{Z}_2,$$
has rank 10.
The intersection pairing $\langle,\rangle$ on $H_2(X,\mathbb{Z})'$ 
is
isomorphic to the quadratic form  
$U\oplus E_8(-1)$,
where 
$$U
= \left( \begin{array}{cc}
0 & 1 \\
1 & 0 \end{array} \right)$$
and 
$$E_8(-1)=  \left( \begin{array}{cccccccc}
 -2&    0 &  1 &   0 &   0 &   0 &   0 & 0\\
    0 &   -2 &   0 &  1 &   0 &   0 &   0 & 0\\
     1 &   0 &   -2 &  1 &   0 &   0 & 0 &  0\\
      0  & 1 &  1 &   -2 &  1 &   0 & 0 & 0\\
      0 &   0 &   0 &  1 &   -2 &  1 & 0&  0\\
      0 &   0&    0 &   0 &  1 &  -2 &  1 & 0\\ 
      0 &   0&    0 &   0 &   0 &  1 &  -2 & 1\\
      0 & 0  & 0 &  0 & 0 & 0 & 1& -2\end{array}\right).$$
is the (negative) Cartan matrix.
A good reference for the classical geometry of $X$ is 
\cite{igor}.

We will view the curve classes of $X$ as lying in  
 $H_2(X,\mathbb{Z})'$.
More precisely,
the integral
\begin{equation}\label{mtor}
\Big\langle \tau_{a_1}(\gamma_1)\cdots \tau_{a_n}(\gamma_n) 
\Big\rangle^X_{g,\beta}
\end{equation}
with $\beta \in H_2(X,\mathbb{Z})'$
is defined to be the {\em sum} of the integrals
associated to the two lifts of $\beta$ to
$H_2(X,{\mathbb Z})$.

Let $Q$ be the Enriques Calabi-Yau 3-fold defined in Section \ref{et45},
\begin{equation}
\label{v23}
Q= \left( K3 \times E\right)/ \langle \sigma \rangle
\end{equation}
Projection to the first factor,
$$\pi_X:Q\rightarrow X,$$
is an elliptic fibration (with all fibers isomorphic to $E$).
The homology mod torsion also projects,
$$\pi_{X*}: H_2(Q,\mathbb{Z})' \rightarrow
H_2(X,\mathbb{Z})'.$$
Since $\pi_X$ has a 0 section $s_0$, there is
decomposition
$$ H_2(Q,\mathbb{Z})'= H_2(X,\mathbb{Z})' \oplus \mathbb{Z} [E], $$
where $[E]$ is class of the fiber of $\pi_X$.
The Gromov-Witten invariant 
$$
N_{g,(\beta,d)} = \langle 1 \rangle^Q_{g, (\beta,d)}
$$
is defined by summation over all lifts of $(\beta,d)\in
H_2(Q,\mathbb{Z})'$ following the convention of \eqref{mtor}.

By stability, the curve class 0 invariants $N_{g,(0,0)}$ are
{\em not} considered in genus 0 or 1. 

\subsubsection{The projection $\pi_{\mathbb{P}^1}$}
Projection to the second factor of \eqref{v23},
$$\pi_{\mathbb{P}^1}: Q \rightarrow \mathbb{P}^1,$$
is a $K3$ fibration with 4 double Enriques fibers.
The homology classes $(\beta,0)$ of $Q$ project to
0 under $\pi_{\mathbb{P}^1*}$.
The invariants $N_{g,(\beta,0)}$
are the {\em fiber class} invariants of $Q$ with respect to
$\pi_{\mathbb{P}^1}$.

\begin{Lemma} \label{fcl} The fiber class invariants of $Q$ are
Hodge integrals in the Gromov-Witten theory of the
Enriques surface $X$,
$$N_{g,(\beta,0)} = 4  \Big\langle
(-1)^{g-1}\lambda_{g-1} \Big\rangle_{g,\beta}^X.$$
\end{Lemma}

\begin{proof}
We use the good degeneration of $Q$ to $R\cup_{K3} R$
discussed in Section \ref{et45}.
By the degeneration formula,
\begin{equation}\label{aaa1}
N_{g,(\beta,0)} = N_{g,(\beta,0)}^{R/K3} + N_{g,(\beta,0)}^{R/K3},
\end{equation}
where the superscript $R/K3$ denotes the relative invariants of
the pair.

By degeneration to the normal cone of $K3\subset R$,
\begin{equation}\label{bbb2}
N_{g,(\beta,0)}^R = N_{g,(\beta,0)}^{R/K3} + 
N_{g,(\tilde{\beta},0)}^{K3\times\mathbb{P}^1/K3},
\end{equation}
where the relative
divisor 
$$K3\subset K3\times\mathbb{P}^1$$ 
is a section. Here, $\tilde{\beta}$
stands for all classes pushing forward to $\beta$. 
By the Leray-Hirsch result of 
\cite{mp} and the triviality of $\mathsf{GW}(K3)$,
 $$N_{g,(\tilde{\beta},0)}^{K3\times \mathbb{P}^1/K3}=0.$$
By localization,
\begin{equation}\label{ccc3}
N_{g,({\beta},0)}^R 
= 2 \Big
\langle (-1)^{g-1}\lambda_{g-1} \Big
\rangle_{g,\beta}^X.
\end{equation}
The Lemma is obtained by combining
\eqref{aaa1}-\eqref{ccc3}.
\end{proof}

Lemma \ref{fcl} is interpreted as a vanishing result for 
the fiber class invariants of $Q$,
\begin{equation}
\label{df12}
N_{0,(\beta,0)} = 0,
\end{equation}
in genus 0.

\subsubsection{Heterotic string}
Klemm and Mari\~no have determined
the fiber class invariants $N_{g,(\beta,0)}$
in terms of modular forms by heterotic string
calculations \cite{klemmm}. Precise formulas
will be discussed in Section \ref{mf}. 
By Lemma \ref{fcl}, the fiber class results of
\cite{klemmm} may be viewed as computing $\lambda_{g-1}$
Hodge integrals in the Gromov-Witten theory of the Enriques surface $X$.

\subsection{Hodge integrals on the Enriques surface for $g\leq 2$}
\label{fcr}
\subsubsection{Overview}
We will study here the fiber class invariants of $Q$
via Gromov-Witten theory
for $g\leq 2$. In genus 0, vanishing has already been obtained
\eqref{df12}. In genus 1,
 $$N_{1,(\beta,0)} = 4 \langle 1 \rangle_{1,\beta}^X$$
by Lemma \ref{fcl}, and
Hodge insertions on $X$ do not arise.
In genus 2,
$$N_{2,(\beta,0)} = -4 \langle \lambda_1 \rangle_{2,\beta}^X,$$
The required genus 1 and 2 invariants of $X$ will be calculated
by a combination of techniques.

\subsubsection{Isotropic classes}
Let $X$ be an Enriques surface presented as an elliptic fibration 
$$f:X \rightarrow \mathbb{P}^1$$
with 12 singular fibers and 2 double fibers.
Let $$F_f\in H_2(X,\mathbb{Z})'$$
denote half the class of the general fiber of $f$.
A class $\beta\in H_2(X,\mathbb{Z})'$ is {\em positive} if
either 
$$\langle F_f,\beta\rangle>0$$
or $\beta$ is a positive multiple
of $F_f$.

The classes of $H_2(X,\mathbb{Z})'$ represented by algebraic curves
are {\em effective}.
Effective classes must be 
positive.

A {\em primitive} class in $H_2(X,\mathbb{Z})'$ is nonzero and
not divisible.
Let $F\in H_2(X,\mathbb{Z})'$ be a positive, primitive, isotropic class.
By the classical theory of Enriques surfaces \cite{igor},
the class $2F$ is the fiber of an elliptic pencil
$$f_F:X\rightarrow \mathbb{P}^1$$
with 12 singular fibers and 2 double fibers.\footnote{We assume 
here $X$ is generic in the moduli of Enriques surfaces.} 
We will compute
the invariants $\langle \lambda_{g-1} \rangle_{g,nF}^X$ via the
good degeneration of $f_F$ to 
$X_1\cup_E E(1)$
discussed in Section \ref{13s}.

There are two cases to consider. If $n$ is odd, 
the degeneration formula yields
\begin{equation*}
\langle \lambda_{g-1}\rangle^X_{g,nF} = 
\langle \lambda_{g-1} \rangle_{g,nF}^{X_1/E}  
\end{equation*}
since $nF$ is not represented on $E(1)$.
If $n$ is even, then
\begin{equation*}
\langle \lambda_{g-1}\rangle^X_{g,nF} = 
\langle \lambda_{g-1} \rangle_{g,nF}^{X_1/E} + 
\langle \lambda_{g-1}\rangle_{g,nF}^{E(1)/E}.
\end{equation*}
There is a good degeneration of an elliptically fibered $K3$ surface to
$E(1)\cup_E E(1)$. Hence, for $n$ even,
\begin{equation*}
\langle \lambda_{g-1}\rangle^{K3}_{g,nF} = 
\langle \lambda_{g-1} \rangle_{g,nF}^{E(1)/E} + 
\langle \lambda_{g-1}\rangle_{g,nF}^{E(1)/E}.
\end{equation*}
By the vanishing of Gromov-Witten invariants of the $K3$ surface,
$$\langle \lambda_{g-1}\rangle_{g,nF}^{E(1)/E}=0.$$
We conclude
$$ \langle \lambda_{g-1}\rangle^X_{g,nF} = 
\langle \lambda_{g-1} \rangle_{g,nF}^{X_1/E} $$
for all $n$

By degeneration to the normal cone of $E\subset X_1$,
\begin{equation}\label{bbx2}
\langle \lambda_{g-1} \rangle_{g,nF}^{X_1} =
\langle \lambda_{g-1} \rangle_{g,nF}^{X_1/E} + 
\langle \lambda_{g-1}\rangle_{g,nF}^{E\times \mathbb{P}^1/E}.
\end{equation}
As in the previous paragraph, the second term on the
right of \eqref{bbx2} is absent in the $n$ odd case.

Localization may be applied to the Gromov-Witten invariants of $X_1$,
\begin{eqnarray*}
\langle \lambda_{g-1} \rangle_{g,nF}^{X_1} 
& = & 2 
\langle (-1)^{g-1}\lambda_{g-1}
\rangle_{g,n}^E \\
& = & 2 \sigma_{-1}(n) \ \delta_{g,1}.
\end{eqnarray*}
Here,
$$\sigma_{-1}(n)= \sum_{i|n} \frac{1}{i}.$$
Similarly, if $n$ is even,
$$\langle 1\rangle_{1,nF}^{E\times \mathbb{P}^1/E} 
= \sigma_{-1}(n/2)\ \delta_{g,1}.$$
Equation \eqref{bbx2} then yields the following results.

\begin{Lemma}\label{v12}
The Gromov-Witten invariants of $X$ in positive
isotropic classes
are determined by:
\begin{align*}
\langle 1 \rangle_{1,nF}^{X} & =  2\sigma_{-1}(n) & &  {\text{$n$ odd}},\\
\langle 1 \rangle_{1,nF}^{X} & =  2\sigma_{-1}(n)- \sigma_{-1}(n/2)
& & {\text{$n$ even}}.
\end{align*}
\end{Lemma}

\begin{Lemma} The fiber class invariants  $N_{g,(\beta,0)}$
of $Q$ vanish for nonzero isotropic classes $\beta$
if $g\geq 2$.
\end{Lemma}

\subsubsection{Genus 1 invariants of the Enriques surface}
We derive a relation which determines all genus 1 invariants of
$X$ using localization equations and the (as yet unproven) Virasoro
constraints for $X$.

Consider the 3-fold $Y= X \times \mathbb{P}^1$. The torus $\com^*$
acts on $Y$ via $\mathbb{P}^1$.  Let 
$$[X_0], [X_\infty] \in H^*_{\com^*}(Y,\mathbb{Q})$$
denote the classes of the fibers of $X$ over the $\com^*$-fixed
points of $\mathbb{P}^1$ with tangent weights $1$ and $-1$
respectively.
Let $\beta \in H_2(X,\mathbb{Z})'$ be a nonzero class. 
Certainly,
\begin{equation}\label{xew}
\Big\langle \tau_1([X_0]^2) \Big \rangle_{2,(\beta,1)}^Y
=
0
\end{equation}
Calculation of \eqref{xew} via localization will yield
a nontrivial equation in the Gromov-Witten theory of $X$.

A straightforward
application of the virtual localization formula \cite{GP} yields
\begin{multline*}
\Big\langle \tau_1([X_0]^2) \Big \rangle_{2,(\beta,1)}^Y
=
2 \Big\langle \tau_2 \Big\rangle_{2,\beta}^X - 4 \Big\langle \lambda_1
\Big\rangle_{2,\beta}^X\\ -   \sum_{\beta_1+\beta_2=\beta}
\Big \langle 1 \Big \rangle_{1,\beta_1}^X 
\Big \langle 1 \Big \rangle_{1,\beta_1}^X  \langle \beta_1,\beta_2\rangle.
\end{multline*}
The second term on the right side can be
evaluated via 
the Hodge removal equation of \cite{FP},
\begin{multline*}
\Big\langle \lambda_1 \Big\rangle_{2,\beta}^X =
\frac{1}{12} \Big\langle \tau_2\Big\rangle_{2,\beta}^X 
+\frac{1}{24} \Big\langle 1\Big\rangle_{1,\beta}^X \langle \beta,\beta \rangle
\\+
\frac{1}{24}
 \sum_{\beta_1+\beta_2=\beta}
\Big \langle 1 \Big \rangle_{1,\beta_1}^X 
\Big \langle 1 \Big \rangle_{1,\beta_2}^X  \langle \beta_1,\beta_2\rangle.
\end{multline*}
The first term on the right side can be evaluated by the Virasoro
constraints,
$$\frac{3}{4} \Big\langle \tau_2 \Big\rangle^X_{2,\beta} =
\frac{1}{8} \Big \langle 1 \Big \rangle_{1,\beta}^X 
\langle \beta,\beta \rangle
+\frac{1}{8} \sum_{\beta_1+\beta_2=\beta}
\Big \langle 1 \Big \rangle_{1,\beta_1}^X 
\Big \langle 1 \Big \rangle_{1,\beta_2}^X  \langle \beta_1,\beta_2\rangle,$$
see \cite{ehx,get,pan}.
Equation \eqref{xew} then 
implies the following result.{\footnote{The derivation
depends upon the conjectural Virasoro constraints for X.}}

\begin{Proposition}\label{k34}
For all nonzero $\beta \in H_2(X,\mathbb{Z})'$,
$$\Big\langle 1 \Big \rangle_{1,\beta}^X \langle \beta,\beta
\rangle= 8 \sum_{\beta_1+\beta_2=\beta}
\Big \langle 1 \Big \rangle_{1,\beta_1}^X 
\Big \langle 1 \Big \rangle_{1,\beta_2}^X  \langle \beta_1,\beta_2\rangle.$$
\end{Proposition}

The sum on the right in Proposition \ref{k34} is taken over
all nontrivial decompositions of $\beta$ into effective classes
$\beta_i$ on $X$. 
Effective decomposition defines a partial ordering on
the set of effective classes which has no infinite descending chains.
Proposition \ref{k34} therefore
uniquely determines the genus 1 Gromov-Witten theory of $X$
in terms of the isotropic invariants
$$\Big\langle 1 \Big \rangle _{1,nF}^X$$
calculated in Lemma \ref{v12}.

\begin{Lemma}\label{dq2}
If $\langle \beta,\beta \rangle<0$, then $\langle 1 \rangle^X_{1,\beta}=0$.
\end{Lemma}

\begin{proof}
Let $\beta$ be an effective curve class on
$X$ satisfying
$\langle \beta,\beta \rangle<0$.
Let $$\beta_1+\beta_2=\beta$$ be a
decomposition into effective classes.
Let $H$ be an ample class on $X$. Since
the intersection form on $H_2(X,\mathbb{R})$
has signature $(1,9)$, the form is negative
definite on $H^\perp$. Let
$$\beta_1=h_1 H + N_1, \ \ \beta_2=h_2 H + N_2$$
where $N_i\in H^\perp$.
Since the classes $\beta_i$ are effective, $h_i>0$.
Since $\beta$ has negative square,
$$
(h_1+h_2)\sqrt{\langle H,H\rangle} < \sqrt{-\langle N_1+N_2,N_1+N_2 \rangle}.$$
By the triangle inequality,
$$h_i\sqrt{\langle H,H \rangle}< \sqrt{-\langle N_i,N_i \rangle}$$
must hold for either $\beta_1$ or $\beta_2$. Hence, 
either 
$\langle \beta_1,\beta_1 \rangle<0$ or
$\langle \beta_2,\beta_2 \rangle<0$.
The Lemma is then obtained by induction on the
partial ordering.
\end{proof}

By Lemma \ref{dq2}, we may rewrite Proposition 1 purely in
terms of the intersection form on $H_2(X,\mathbb{Z})'$ ---
without regard to effectivity.

\vspace{10pt}
\noindent{\bf Proposition $\mathbf{1'}$.}
{\em For all nonzero $\beta \in H_2(X,\mathbb{Z})'$,
$$\Big\langle 1 \Big \rangle_{1,\beta}^X \langle \beta,\beta
\rangle= 8 \sum_{\beta_1+\beta_2=\beta}
\Big \langle 1 \Big \rangle_{1,\beta_1}^X 
\Big \langle 1 \Big \rangle_{1,\beta_2}^X  \langle \beta_1,\beta_2\rangle$$
where the sum is over decompositions into
positive classes satisfying $\langle \beta_i,\beta_i\rangle\geq 0$}.
\vspace{10pt}
\label{rt3}
\subsubsection{Genus 2 fiber classes}
The Hodge removal and Virasoro equations of Section \ref{rt3}
yield
$$
\Big\langle \lambda_1 \Big\rangle_{2,\beta}^X = 
\frac{1}{18} \Big\langle 1\Big\rangle_{1,\beta}^X \langle \beta,\beta \rangle
+\frac{1}{18}
 \sum_{\beta_1+\beta_2=\beta}
\Big \langle 1 \Big \rangle_{1,\beta_1}^X 
\Big \langle 1 \Big \rangle_{1,\beta_2}^X  \langle \beta_1,\beta_2\rangle.
$$
After applying Proposition 1,
$$\Big\langle \lambda_1 \Big\rangle_{2,\beta}^X = 
\frac{1}{16} \Big\langle 1\Big\rangle_{1,\beta}^X 
\langle \beta,\beta \rangle.$$
In terms of the invariants of $Q$,
\begin{equation}\label{g34}
N_{2,(\beta,0)} = -\frac{1}{16} 
N_{1,(\beta,0)} \langle \beta,\beta \rangle.
\end{equation}

The Eisenstein series $E_{2n}$ is the modular form
defined by the equation 
$$-\frac{B_{2n}}{4n} E_{2n}(q) = -\frac{B_{2n}}{4n} + \sum_{k\geq 1} 
\sigma_n(k) q^{k},$$
where $\sigma_n(k)$ is the sum of the $n^{th}$ powers of
the divisors of $k$,
$$\sigma_n(k) = \sum_{i|k} i^n.$$
The natural regularization of $\sigma_1(0)$ is
$$\sigma_1(0)=-\frac{B_2}{4}=- \frac{1}{24}.$$
We may rewrite \eqref{g34} as
\begin{equation}\label{f56}
N_{2,(\beta,0)} = \frac{3}{2}\sigma_1(0)  
N_{1,(\beta,0)} \langle \beta,\beta \rangle.
\end{equation}

\subsubsection{Modular forms}\label{mf}

Let $v_1, \ldots v_{10} \in H_2(X,\mathbb{Z})'$ be a basis
with 
$$v_1,v_2\in U, \ \ v_3, \ldots, v_{10}\in E_8(-1)$$
with respect to an identification 
$$H_2(X,\mathbb{Z})' \stackrel{\sim}{=} U \oplus E_8(-1).$$
Let $$v(t)= \sum_{i=1}^nt_i v_i$$
be coordinates in the basis.
Since $v_1$ is a primitive isotropic class,
 positivity can be defined by intersection with
$v_1$. A vector 
$$\beta= \sum_{i=1}^n b_i v_i$$ is 
positive if $b_2>0$ {\em or} if $b_2=0$ and $\beta$ 
is a positive
multiple of $v_1$.

The fiber class potential function of the Enriques
Calabi-Yau 3-fold is defined by
$$F_g(t) = \sum_{\beta>0} N_{g,(\beta,0)} e^{-\langle \beta, v(t) \rangle}.$$
The heterotic string evaluation of $F_g$
by Klemm and Mari\~no  
is 
\begin{equation*}
F_g(t) = \sum_{\beta >0} c_g\big(\langle \beta, \beta \rangle\big)
\Big( 2^{3-2g} {\text{Li}}_{3-2g}( e^{-\langle \beta, v(t) \rangle}) -
{\text{Li}}_{3-2g}( e^{-2\langle \beta, v(t) \rangle}) \Big),
\end{equation*}
for $g>0$.

The terms on the right side are explained as follows.
The coefficients $c_g(n)$ are defined by
$$\sum_{n} c_g(n)q^n =
-\frac{2}{q} \prod_{n=1}^\infty (1-q^{2n})^{-12} \cdot {\mathcal P}_g(q).$$
Here, ${\mathcal P}_g(q)$ is a quasimodular form,
$${\mathcal P}_g(q)= {\mathcal S}_g\left( x_k= \frac{|B_{2k}|}{(2k)!}
E_{2k}(q)\right),$$
with the polynomial $\mathcal{S}_g(x_1,\ldots,x_g)$ defined by
$$
\sum_{g=0}^\infty \mathcal{S}_g(x_1, \ldots, x_g) z^g=
\exp\left[ \sum_{g=1}^\infty x_g z^g \right].$$
For example,
$${\mathcal P}_1 = \frac{1}{12} E_2, \ \  {\mathcal P}_2 =\frac{1}{1440}
(5 E_2^2 + E_4).$$
The polylogarithm $\text{Li}_k$ is defined by
$$\text{Li}_k(x) = \sum_{n=1}^\infty
\frac{x^n}{n^k}.$$

Proposition $1'$ should yield a coefficient relation
for the various modular series in $F_1(t)$.
We have not yet fully checked the compatibility.
Many parallel properties hold. For example,
equation \eqref{f56} is valid for the Klemm-Mari\~no
formula, see Section 4.1 of \cite{klemmm}.

\subsubsection{Donaldson-Thomas theory}
By the GW/DT correspondence, we may instead study the Donaldson-Thomas
theory of $Q$.
The fiber
Donaldson-Thomas theory of $Q$  has a
reduction to the classical cohomology of the 
Hilbert scheme of points of $X$.
We expect the resulting vertex algebra calculations to be very
closely related to the heterotic string results of \cite{klemmm}.

\subsection{Vanishing in genus 0 and 1}
We now consider the Gromov-Witten theory of $Q$ for all
classes $(\beta,d) \in H_2(Q,\mathbb{Z})'$.

\begin{Proposition}
$N_{0,(\beta,d)}=0$.\label{gte}
\end{Proposition}

\begin{proof}
If $d=0$, the vanishing \eqref{df12} has already been been
obtained.
We assume $d>0$.

We use the good degeneration of $Q$ to $R\cup_{K3} R$
discussed in Section \ref{et45}.
By the degeneration formula,
\begin{equation}\label{aaaf1}
N_{0,(\beta,d)} = \sum_\eta \sum_{\beta_1+\beta_2=\beta}
\Big \langle \ 1 \ \Big | 
 \eta \Big \rangle_{g_1,(\beta_1,d)}^{\bullet R/K3}
\Big \langle \eta^\vee \Big | 
 \ 1\ \Big \rangle_{g_2,(\beta_2,d)}^{\bullet R/K3}
\end{equation}
where $\eta^\vee$ is the partition with cohomology weights
Poincar\'e dual to the cohomology weights of $\eta$.

The left side of \eqref{aaaf1} is a connected invariant.
The superscript $\bullet$ on the right side of
\eqref{aaaf1} denotes disconnected invariants. After 
 connecting the
disconnected domains via the gluing conditions specified
by $(\eta,\eta^\vee)$, a  
connected domain must be obtained.
Also, the genus condition
\begin{equation}\label{xdf}
g_1 + g_2+ \ell(\eta)-1 = 0
\end{equation}
must be satisfied.

The relative invariant $\Big \langle \ 1 \ \Big | 
 \eta \Big \rangle_{0,(\beta_1,d)}^{\bullet R/K3}$ can be calculated
from the absolute invariants of $R$ and $K3$ by Theorem 2 
\cite{mp}. Since the Gromov-Witten invariants of $R$ and $K3$ 
vanish in genus 0 for classes with nonzero push-forwards to
$X$, 
we conclude
$$\Big \langle \ 1 \ \Big | 
 \eta \Big \rangle_{0,(\beta_1,d)}^{\bullet R/K3} = 0$$
unless $\beta_1=0$.
Hence, $N_{0,(\beta,d)}$ vanishes if $\beta\neq 0$.

If $\beta=0$, the degeneration formula \eqref{aaaf1}
yields a vanishing for a more subtle reason.
The relative divisor $K3\subset R$ intersects every
fiber of 
$$\pi: R \rightarrow X$$
in exactly 2 points.
Hence, $\ell(\eta)$ must be at least twice the number of connected
components of the domain in a nonvanishing invariant
$\Big \langle \ 1 \ \Big | 
 \eta \Big \rangle_{0,(\beta_1,d)}^{\bullet R/K3}$.
Since each connected domain component must be of genus 0,
the genus condition \eqref{xdf} can never be satisfied. 
\end{proof}

In genus 1, a similar vanishing result holds for classes which
are a proper mix of fiber and base curves.

\begin{Proposition}
For $\beta\neq 0$ and $d\neq 0$,
$N_{1,(\beta,d)}=0.$
\end{Proposition}

\begin{proof}
Consider the degeneration formula
\begin{equation}\label{aaax1}
N_{1,(\beta,d)} = \sum_\eta \sum_{\beta_1+\beta_2=\beta}
\Big \langle \ 1 \ \Big | 
 \eta \Big \rangle_{g_1,(\beta_1,d)}^{\bullet R/K3}
\Big \langle \eta^\vee \Big | 
 \ 1\ \Big \rangle_{g_2,(\beta_2,d)}^{\bullet R/K3}.
\end{equation}
If the connected components of the relative invariants
all have domain genus 0, then $\beta=0$ as in the proof
of Proposition \ref{gte}. If 
a connected component of genus 1 occurs, then the
genus condition
$$
g_1 + g_2+ \ell(\eta)-1 = 1
$$
can not be satisfied unless $d=0$.
\end{proof}

\begin{Lemma} For $d>0$,
$N_{1,(0,d)}= 12\sigma_{-1}(d)$.
\end{Lemma}
\begin{proof}
The Lemma may be obtained by an elementary evaluation of the
degeneration formula \eqref{aaax1} or from the Gromov-Witten
theory of the fibration
$\pi_X: Q \rightarrow X.$
For the latter derivation,
\begin{eqnarray*}
N_{1,(0,d)}& =&  \int_{[\overline{M}_{1}(Q,d[E])]^{vir_\pi}} 
c_2(\mathbb{E}^\vee \otimes T_X) \\
& = & 12 \int_{[\overline{M}_{1}(E,d[E])]^{vir}} 1\\
& = & 12 \sigma_{-1}(d). 
\end{eqnarray*}
Here, $vir_\pi$ is the relative virtual class of the morphism $\pi$,
and $\mathbb{E}$ is the Hodge bundle.
\end{proof}

Since the absolute Gromov-Witten theory of a target
 elliptic curve $E$ {\em with trivial integrand} vanishing in genus $g\geq 2$,
the above proof also yields the following result.

\begin{Proposition}\label{drt}
$N_{g,(0,d)}=0$ for $g\geq 2$ and $d\geq 0$.
\end{Proposition}

\subsection{Holomorphic anomaly in genus 2}
\subsubsection{Overview}
Our last topic is the Gromov-Witten theory of  $Q$ in genus 2.
The holomorphic anomaly equation arises naturally in our geometric
study.
\subsubsection{Absolute invariants}
We will require the explicit evaluations of several 
genus 1 invariants of  $R$.
Let $$\iota: K3 \rightarrow R$$
denote the inclusion of the relative divisor, and let
$$\pi:K3\rightarrow X$$
denote the induced projection.
As in Section \ref{fcr}, the small brackets $\langle,\rangle$
denote the intersection pairing on $X$.

\begin{Lemma} For $\gamma\in H_2(K3,\mathbb{Q})$,
$$\Big \langle \tau_{2d-1}(\iota_*(\gamma)) \Big \rangle_{1,(\beta,d)}^{R}
= \frac{2d}{(d!)^2}
\Big \langle 1 \Big \rangle_{1,\beta}^X  \langle \beta,\pi_*(\gamma) 
\rangle.
$$\label{d5}
\end{Lemma}

\begin{proof}
By parallel localization arguments, we find
$$
\Big \langle \tau_{2d-1}(\iota_*(\gamma)) \Big \rangle_{1,(\beta,d)}^{R}
=
\Big \langle \tau_{2d-1}(p) \Big \rangle_{1,(1,d)}^{E\times \mathbb{P}^1}
\Big \langle 1 \Big \rangle_{1,\beta}^X  \langle \beta,\pi_*(\gamma) 
\rangle
$$
where $p$ is a point class on the surface $E\times \mathbb{P}^1$.
By degenerating the elliptic factor $E$ of the target
 to a nodal rational curve,
$$\Big \langle \tau_{2d-1}(p) \Big \rangle_{1,(1,d)}^{E\times \mathbb{P}^1}=
2 \Big \langle (1,p') \ \Big| \ \tau_{2d-1}(p) \ \Big|
\ (1,1) \Big \rangle_{0,(1,d)}
^{\mathbb{P}^1 \times \mathbb{P}^1/ \mathbb{P}^1_0\cup
\mathbb{P}^1_\infty}$$
where $p'$ is a point class of the relative divisor $\mathbb{P}^1_0$.
The exchange of relative for absolute insertions takes a simple
form here,
$$\Big \langle (1,p') \ \Big| \ \tau_{2d-1}(p) \ \Big|
\ (1,1) \Big \rangle_{0,(1,d)}
^{\mathbb{P}^1 \times \mathbb{P}^1/ \mathbb{P}^1_0\cup
\mathbb{P}^1_\infty} =
\Big \langle \tau_0(p') \tau_{2d-1}(p) \tau_0([\mathbb{P}^1_\infty])
\Big \rangle_{0,(1,d)}
^{\mathbb{P}^1 \times \mathbb{P}^1}.$$
We apply topological recursion relations to the result,
\begin{multline*}
\Big \langle \tau_0(p') \tau_{2d-1}(p) \tau_0([\mathbb{P}^1_\infty])
\Big \rangle^{\mathbb{P}^1 \times
\mathbb{P}^1}_{0,(1,d)} =\\
\Big \langle \tau_{2d-2}(p) \tau_0([1,0])
\Big \rangle^{\mathbb{P}^1 \times
\mathbb{P}^1}_{0,(0,d)}
\Big \langle \tau_0([0,1]) \tau_0(p')  \tau_0([0,1])
\Big \rangle^{\mathbb{P}^1 \times
\mathbb{P}^1}_{0,(1,0)}
\end{multline*}
where $[1,0],[0,1]\in H^2(\mathbb{P}^1 \times
\mathbb{P}^1, \mathbb{Q})$ denote the classes of the 2 rulings.
The evaluations
\begin{eqnarray*}
\Big \langle \tau_{2d-2}(p) \tau_0([1,0])
\Big \rangle^{\mathbb{P}^1 \times
\mathbb{P}^1}_{0,(0,d)} & = &
d \Big \langle \tau_{2d-2}(p) 
\Big \rangle^{\mathbb{P}^1}_{0,d} \\
& = & \frac{d}{(d!)^2}
\end{eqnarray*}
and
$$\langle \tau_0([0,1] \tau_0(p')  \tau_0([0,1])
\Big \rangle^{\mathbb{P}^1 \times
\mathbb{P}^1}_{0,(1,0)}=1$$
complete the derivation.

The well-known evaluation of the genus 0
descendent $\langle \tau_{2d-2}(p) \rangle^{\mathbb{P}^1}_{0,d}$ used above
can be
found in \cite{pgiv}.
\end{proof}

Almost identical arguments yield the evaluations of 
the following 2-point invariants.

\begin{Lemma} \label{d6}
Let $a_1=2m_1+1$ and $a_2=2m_2+1$ be odd integers
satisfying
$a_1+a_2=2d$.
For $\gamma_1,\gamma_2\in H_2(K3,\mathbb{Q})$,
\begin{multline*}\Big \langle \tau_{a_1-1}(\iota_*(\gamma_1)) 
\tau_{a_2-1}(\iota_*(\gamma_2)) 
\Big \rangle_{1,(\beta,d)}^{R}
=\\ \frac{2}{(m_1!)^2(m_2!)^2}
\Big \langle 1 \Big \rangle_{1,\beta}^X  \langle \beta,\pi_*(\gamma_1) 
\rangle\langle \beta,\pi_*(\gamma_2) 
\rangle
.
\end{multline*}
\end{Lemma}

\begin{Lemma} \label{d7}
Let $a_1=2m_1$ and $a_2=2m_2$ be even integers
satisfying
$a_1+a_2=2d$.
For $\gamma_1,\gamma_2\in H_2(K3,\mathbb{Q})$,
\begin{multline*}\Big \langle \tau_{a_1-1}(\iota_*(\gamma_1)) 
\tau_{a_2-1}(\iota_*(\gamma_2)) 
\Big \rangle_{1,(\beta,d)}^{R}
=\\ \frac{2mn}{(m_1!)^2(m_2!)^2}
\Big \langle 1 \Big \rangle_{1,\beta}^X  \langle \beta,\pi_*(\gamma_1) 
\rangle\langle \beta,\pi_*(\gamma_2) 
\rangle
.
\end{multline*}
\end{Lemma}

\subsubsection{Relative invariants}
The degeneration to the normal cone of $K3\subset R$
can be applied to determine relative invariants from 
absolute invariants, see Theorem 2 \cite{mp}.

\begin{Lemma} For $\gamma\in H_2(K3,\mathbb{Q})$,
$$\Big\langle 1 \ \Big |\ (2d,\gamma) \ \Big \rangle_{1,(\beta,d)}^{R/K3} =
2 
\Big \langle 1 \Big \rangle_{1,\beta}^X  \langle \beta,\pi_*(\gamma) 
\rangle.$$ \label{w1}
\end{Lemma}

\begin{proof}
Let $I_d$ denote the relative invariants to be determined,
$$I_d=\Big\langle 1 \ \Big |\ (2d,\gamma) \ 
\Big \rangle_{1,(\beta,d)}^{R/K3}.$$
Degeneration to the normal cone of $K3\subset R$ yields
\begin{multline*}
\Big \langle \tau_{2d-1}(\iota_*(\gamma)) \Big \rangle_{1,(\beta,d)}^{R} = 
I_d \ 2d \ \Big \langle (2d) \ \Big | \ \tau_{2d-1}(p) \Big 
\rangle_{0,2d}^{\mathbb{P}^1} \\ +
\sum_{r=1}^{d-1} I_r \ 2r \ \Big \langle (2r),(1)^{d-r} \ \Big | \ \tau_{2d-1}(p) \Big 
\rangle_{0,d+r}^{\mathbb{P}^1}.
\end{multline*}
The $K3\times \mathbb{P}^1/K3$ side has been written in terms
of the relative Gromov-Witten theory of a vertical $\mathbb{P}^1$
since $\mathsf{GW}(K3)$ is trivial.
The coefficients
$$\Big \langle (2d) \ \Big | \ \tau_{2d-1}(p) \Big 
\rangle_{0,2d}^{\mathbb{P}^1} =\frac{1}{(2d)!},$$
$$\Big \langle (2r),(1)^{d-r} \ \Big | \ \tau_{2d-1}(p) \Big 
\rangle_{0,d+r}^{\mathbb{P}^1} = \frac{1}{(d+r)!(d-r)!}$$
are easily evaluated by completed cycles \cite{OP1}.
We find
$$I_d = 2 
\Big \langle 1 \Big \rangle_{1,\beta}^X  \langle \beta,\pi_*(\gamma) 
\rangle$$
is the unique solution to the recursion
\begin{equation*}
\frac{2d}{(d!)^2}
\Big \langle 1 \Big \rangle_{1,\beta}^X  \langle \beta,\pi_*(\gamma) 
\rangle = 
\frac{I_d}{(2d-1)!} + \sum_{r=1}^{d-1} \frac{2rI_r}{(d+r)!(d-r)!} 
\end{equation*}
obtained from Lemma \ref{d5}.
\end{proof}

A parallel (though more complicated) 
derivation from the evaluations of Lemmas \ref{d6} and \ref{d7}
yields a second result.

\begin{Lemma} For $\gamma_1,\gamma_2 \in H_2(K3,\mathbb{Q})$,
$$\Big\langle 1 \ \Big |\ (d,\gamma_1),(d,\gamma_2) 
\ \Big \rangle_{1,(\beta,d)}^{R/K3} =
2 
\Big \langle 1 \Big \rangle_{1,\beta}^X  \langle \beta,\pi_*(\gamma_1) 
\rangle
\langle \beta,\pi_*(\gamma_2) 
\rangle
.$$\label{w2}
\end{Lemma}

Lemmas \ref{w1} and \ref{w2} have a remarkable property --- the
relative evaluations are independent of $d$.

\subsubsection{Holomorphic anomaly}
Since the invariant
$N_{2,(\beta,d)}$ has been determined by\  Lemma \ref{drt} if
$\beta=0$ and by 
\eqref{f56}
if $d=0$, we assume $\beta \neq 0$ and $d>0$.

We calculate $N_{2,(\beta,d)}$ via the degeneration of
$Q$ to $R\cup_{K3} R$,
\begin{equation}\label{xxxxy}
N_{2,(\beta,d)} = \sum_\eta \sum_{\beta_1+\beta_2=\beta}
\Big \langle \ 1 \ \Big | 
 \eta \Big \rangle_{g_1,(\beta_1,d)}^{\bullet R/K3}
\Big \langle \eta^\vee \Big | 
 \ 1\ \Big \rangle_{g_2,(\beta_2,d)}^{\bullet R/K3}.
\end{equation}
If a connected
component 
of genus 2 occurs on the right side of degeneration formula, then all other
components must be genus 0 vertical classes. Since each genus 0
vertical class must intersect the relative divisor at least
twice, the genus condition 
$$g_1 + g_2+ \ell(\eta)-1 = 2$$
can not be satisfied.

Since $\beta \neq 0$, a connected component of genus 1 must occur
on the right side of degeneration formula
\eqref{xxxxy}.
There are exactly two possibilities:
\begin{enumerate}
\item[(i)] the degeneration graph has a single genus 1 component
           with a self node, 
\item[(ii)] the degeneration graph has two genus 1 components.
\end{enumerate}
Genus reduction is the hallmark of the holomorphic anomaly 
equation.

Consider first the degeneration terms of type (i). The elliptic
component may occur on either side of $R\cup_{K3} R$.
All other components must be genus 0 vertical curves
fully ramified at the intersection points with the relative
divisor.
Once the side of the elliptic component is specified,
the geometric configurations are easily seen to be in bijective
correspondence with divisors of $d$.
The term corresponding to the divisor $r$ contributes
\begin{eqnarray*}
\sum_{i}
\frac{r}{2} 
\Big\langle 1 \ \Big |\ (r,\gamma_i),(r,\gamma^\vee_i) 
\ \Big \rangle_{1,(\beta,d)}^{R/K3}&  = &
\sum_i r
\Big \langle 1 \Big \rangle_{1,\beta}^X  \langle \beta,\pi_*(\gamma_i) 
\rangle
\langle \beta,\pi_*(\gamma^\vee_i) 
\rangle \\
&= & 2r \Big \langle 1 \Big \rangle_{1,\beta}^X \langle \beta,\beta \rangle,
\end{eqnarray*}
where the sum is over a basis $\{\gamma_1, \ldots, \gamma_{22}\}$
 of $H_2(K3,\mathbb{Q})$. Lemma \ref{w2}
is used for the first equality.

The full type (i) contribution to the degeneration formula \eqref{xxxxy} is
\begin{equation}\label{h11}
4 \sigma_1(d) \ \Big \langle 1 
\Big \rangle_{1,\beta}^X \langle \beta,\beta \rangle,
\end{equation}
counting both sides for the elliptic component.

The degeneration terms of type (ii) have a similar treatment.
Again, a bijective correspondence with divisors of $d$ is found.
The full type (ii) contribution to \eqref{xxxxy} is
\begin{equation}\label{h22}
\sum_{\beta_1+\beta_2=\beta} 16 \sigma_1(d)\ \Big \langle 1 \Big 
\rangle_{1,\beta_1}^X
\Big \langle 1 \Big \rangle_{1,\beta_2}^X \langle \beta_1,\beta_2 \rangle
\end{equation}
using the evaluation of Lemma \ref{w1}.

Summing the contributions and writing the result in terms of
the fiber class Gromov-Witten invariants of $Q$ by Lemma \ref{fcl}
yields the following result.

\begin{Theorem} For $d>0$, \label{ve2}
\begin{equation*}
N_{2,(\beta,d)} = \sigma_1(d) \Big( N_{1,(\beta,0)} \langle \beta,\beta
\rangle +
\sum_{\beta_1+\beta_2=\beta} N_{1,(\beta_1,0)} N_{1,(\beta_2,0)} 
\langle \beta_1,\beta_2\rangle 
\Big).
\end{equation*}
\end{Theorem}

Theorem 3 may be interpreted as the holomorphic anomaly equation
in genus 2 for the Enriques Calabi-Yau 3-fold.{\footnote{Our
derivation of Theorem 3 does {\em not} depend upon the
Virasoro constraints for the Enriques surface}} A discussion
can be found in Sections 6.2.2 - 6.2.4 of \cite{klemmm}.
In fact, Theorem 3 is used in \cite{klemmm} to fix the
holomorphic ambiguity.

We may rewrite Theorem \ref{ve2} using the fiber class
results in genus 1 and 2 of Section \ref{fcr}. By Proposition 1,
$$
N_{2,(\beta,d)}  =  \frac{3}{2} \sigma_1(d) N_{1,(\beta,0)} \langle
\beta, \beta \rangle
$$
for $d>0$.
By \eqref{f56},
$$N_{2,(\beta,0)} = \frac{3}{2} \sigma_1(0) N_{1,(\beta,0)} \langle
\beta,\beta \rangle.$$

\begin{Corollary} We have
\begin{eqnarray*}
\sum_{d\geq 0} N_{2,(\beta,d)} q^d & =&  -\frac{1}{16}
E_2(q) N_{1,(\beta,0)} \langle \beta,\beta\rangle \\
& = & E_2(q) N_{2(\beta,0)}.
\end{eqnarray*}
\end{Corollary}

We have calculated the Gromov-Witten theory of $Q$
in genus $g\leq 2$. We expect the Gromov-Witten theory
is exactly solvable in all genera.

The Enriques Calabi-Yau 3-fold may be the most tractable
compact Calabi-Yau with nontrivial Gromov-Witten theory. 
Certainly the higher genus study of the quintic 3-fold
in $\mathbb{P}^4$ appears more difficult, see \cite{Ga,LZ,mp}.

\vspace{+10 pt}
\noindent
Department of Mathematics \\
Princeton University \\
Princeton, NJ 08544, USA\\
dmaulik@math.princeton.edu \\

\vspace{+10 pt}
\noindent
Department of Mathematics\\
Princeton University\\
Princeton, NJ 08544, USA\\
rahulp@math.princeton.edu

\end{document}